# The Role of Professional Societies in STEM Diversity

Vernon R. Morris and Talitha M. Washington, Howard University, Washington, DC
Email: vmorris@howard.edu, talitha.washington@howard.edu

**Abstract**

The overall percentages of African American scientists indicate underrepresentation in most science, technology, engineering, and mathematics (STEM) disciplines and the percentages appear to be declining over the last three decades [NSF, 2017]. Despite investments in diversity programs, the observable impact on STEM leadership and the demographics of the science and technology workforce remains quite small. This presentation will highlight some of the challenges and barriers that many students and professionals who seek to pursue careers in these fields face, and the role of professional societies in either exacerbating the perpetuation of monocultures in the various STEM disciplines or proactively working to eliminate barriers and discrimination. We will present and provide clarity on three common myths that are often articulated in discussions of STEM diversity. We will share insights on how professional societies can directly impact the broadening of participation as well as the persistence of racial groups in the STEM fields and hence, strengthen and sustain the Nation's future workforce.

**Background**

Both authors have been active members of large majority- and minority-based professional scientific organizations for the past twenty years. Within the majority professional societies, we have served on numerous ad hoc organizations and committees committed to addressing access and inclusion issues that are common to both majority and minority organizations. The focus of this article is on racial diversity; in particular, the participation of African Americans in STEM and potential ways scientific organizations can aid in the diversification of STEM.

**What is the current picture?**

This paper focuses on the Geosciences and Mathematics disciplines given that they are characterized by severe underrepresentation. A recent New York Times article underscored the ineffectiveness of Affirmative Action programs in making a sizeable or sustainable impact on STEM diversity at the top colleges in the Nation over the past 35 years (Ashkenas et al. 2017). Moreover, the article shows trends for participation by African Americans and Hispanics were negative while the trends for Caucasians and Asians were both positive over the same period. Given the demographic patterns among the current and future student populations, this is not just an alarming observation, it is an issue of National Security. Home-grown expertise in STEM undergirds the need for continued advances in cybersecurity, food security, and environmental intelligence National Defense and Homeland Security.

Figure 1 shows the 2006 participation within science and engineering occupations for three of the major ethnic groups considered to be perennially underrepresented in STEM fields. While these data are a decade old, the trends in participation over the past decade have remained negative (NSF 2017). Thus, we can treat these percentages as an effective upper bound to the actual numbers. The total participation of underrepresented minorities in Geosciences and Computer and Mathematical Sciences are 6% and 8.9%, respectively. We note that both of these fields have been identified as the fastest growing jobs in the sciences (Eos 2017). According to the 2016 census population estimates, the percentage of Black or African American is 13.3%, Hispanic or Latino is 17.8%, and the percentage of Native American, Alaska Nave, Native Hawaiian and Other Pacific Islander is 1.5% (US Census Report 2015). If parity in the fields were achieved, then the percentages in the fields would be a little closer to 32.6%, not 6% or 8.9%. Comparing these percentages for underrepresented groups, the problem of underrepresentation in STEM is clear.

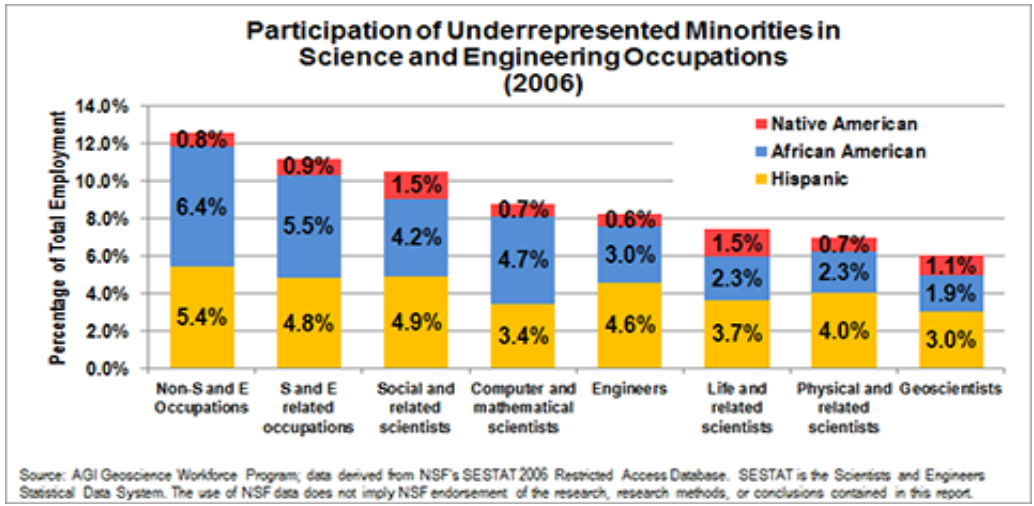

Figure 1: Participation of Underrepresented Minorities in Science
and Engineering Occupations (ACI Workforce Report 2006)

Despite the fact that the total number of blacks earning bachelor's degrees has increased by 54%, the fraction relative to the total number of those earning a STEM bachelor's degree has only increased from 10% to 11%. Meanwhile, the percentage of STEM degrees earned by blacks appears to be decreasing. That is even though there are more blacks earning bachelor's degrees, fewer of them are earning bachelor's in STEM disciplines (NSF 2017).

Across racial and ethnic groups the recent decadal trend also indicates larger shares of undergraduate degrees and certificates being conferred to female students than to male students (NSF 2017). We recognize that the trend that we see at the postsecondary level impacted by the rates of minority male incarceration and rates of high school attrition. We remain cognizant of the entire educational experience of minorities in STEM and at the moment focus on the collegiate and professional levels.

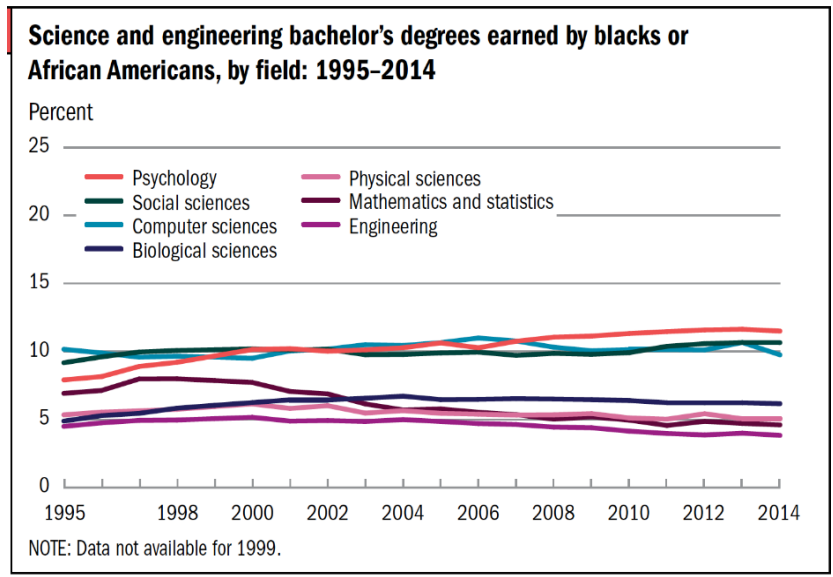

Figure 2: Science and engineering bachelor's degrees earned by Blacks or African-Americans (NSF 2017)

In Figure 2, the trendline for black or African American bachelor's degrees in mathematics and statistics is decreasing over the time interval from 1995 to 2014, from about 7% to less than 5%. As the geosciences is in the physical sciences, Figure X indicates that the percentage of black or African Americans earning bachelor's degrees in this field is relatively flat, around 5%. Note that these percentages are significantly less than the percentages of black or African Americans in the US population in 2016, which is 13.3% (NSF, 2017).

Nearly 30% of all black S&E students receiving their doctoral degrees from US universities earned their bachelor's degrees at an HBCU (NSF 2017) even though only 9% of blacks attended HBCUs in 2015 (Pew Center 2017). HBCUs are even more overrepresented as the "incubators" of future doctoral degree recipients in mathematics and geosciences (e.g. Earth, Atmospheric and Ocean Sciences). As an example, in 2014, over 40% of the African American students in atmospheric sciences were either educated at an HBCU or earned their terminal degree at an HBCU. As shown in Figure 3, many graduate STEM programs are heavily populated by blacks who attend HBCUs. However, over the years 2004 to 2014, the percentage of black or African Americans earning S&E degrees from HBCUs is declining (see Figure 4) which raises great concern about the future of blacks in STEM.

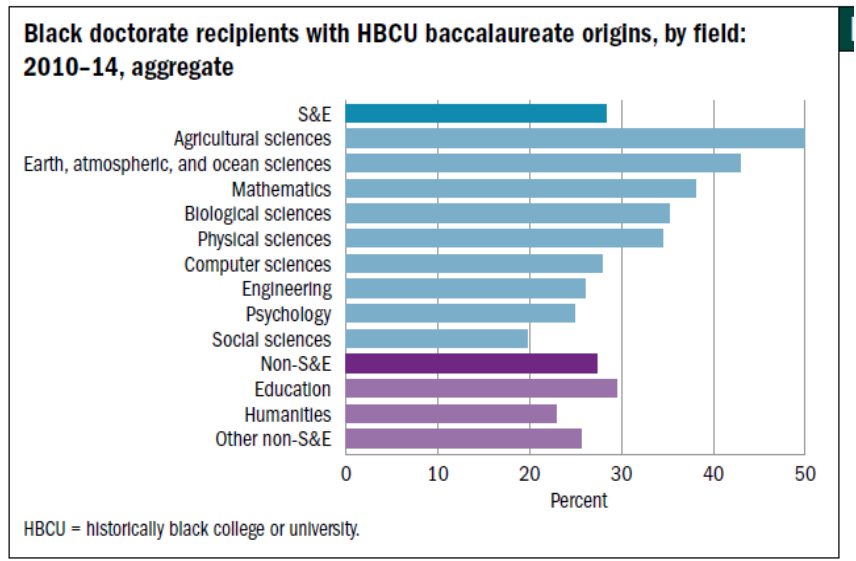

Figure 3: Black doctorate recipients who graduated from HBCUs (NSF 2017)

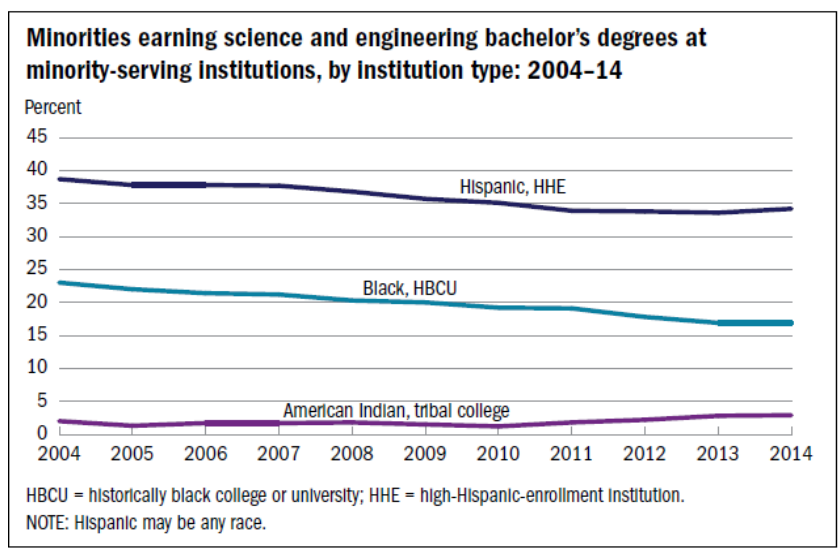

Figure 4: Minorities earning science and engineering bachelor's degrees at minority serving institutions (NSF 2017)

Professional societies often play an intermediary role between students seeking careers and opportunities within the discipline. For the geosciences and mathematical sciences; two disciplines that are among the least diverse by NSF indicators, there are longstanding cultures that have developed within federal agencies, elite academic programs, and amongst private sector stakeholders that tend towards homogenization rather than diversity. We recognize that there are a plethora of ways that diversity in STEM access and inclusion can be improved. This paper focuses on one key element, the scientific professional societies, in an attempt to significantly add to the ongoing holistic strategies that seek to improve retention and growth of diversity in various STEM fields.

**Personal stories**

One co-author, Morris, can still remember experiences during the first professional meetings he attended over 25 years ago. He recalls walking anxiously and excitedly through the crowded venues encountering professionals whose papers he had read voraciously, hoping to steal a few moments of their time and garner a prized opinion or word of encouragement. More often than not, they would avert their eyes (even when talking to him), refrain from shaking hands on introductions, or ignore him completely. Often, he has traveled through the crowded meeting halls and instead of being greeted by his peers, they withdraw from him like "the parting of the Red Sea". When these experiences of how our Caucasian attendees respond to visually ethnic minorities in these crowded professional or social spaces are related to African American and Latino/a colleagues of all genders, it is often received with affirming and empathizing nods because so many find it to be a common encounter. In contrast, he recalls attending meetings of minority professional societies such as the SACNAS Conference where a random lone person would immediately attract other attendees simply seeking to make a scientific connection. Why this same spirit of acceptance cannot be achieved at all meetings is baffling.

The other co-author, Washington, recalls chatting at a breakout table at a conference with her colleague, Dr. Ronald Mickens, Calloway Professor of Physics at Clark Atlanta University. While pontificating about the discretization of dynamical systems, an elderly white woman rolls up in an electric chair. We casually looked over and she asked Professor Mickens, "Do you work here at the hotel?" He replied, "Why do you ask?" Looking around, the lady said, "I need to know where I can put my chair." We all looked at each other baffled since we were surrounded by other scientists and we had our conference badges proudly on display. Smiling, Mickens turns to the lady and excitingly waves his hand around the room to say, "Any where you want." The lady was pleased and rolled away. We began to try to understand why she thought this man, clearly a conference attendee, a scientist who has has published over 300 scientific research papers and over 17 books, could be confused for being a hotel worker. Unfortunately, this is the experience of many professionals of color. We are victims of cognitive dissonance - we cannot be scientists and engineers even in spaces occupied by scientists and engineers. It is easier to assume that we are the "help". In a moment, the sense of belonging and scientific confidence can be shaken down to its core. These assumptions clearly convey who belongs and who does not belong. Why belonging to the scientific academy is assumed to belong to whites only is baffling.

Both of the anecdotes above describe how the monoculture treats or responds to individuals they conceive of as "other". While individuals of past generations or significant stature can shrug off these incidents, it may be much harder for younger colleagues and certainly for post-secondary students. What often results is the dilemma of the imposter syndrome which is a concept that describes the sense of alienation that minorities feel when they are subjected to the types of microaggressions described above. Even high-performing and competent students and professionals can suffer from imposter syndrome. That is, they can begin to suspect that their achievements are the result of luck or accident rather than their own talent and actions and see "others" as inherently more worthy.

Social psychologist Amy Cuddy presented a TED talk in which she explains how nonverbal expressions of power and dominance influence how we think and feel about ourselves as well as how others see us (Fiske et al. 2002). She conveys that it can be very difficult for people to see themselves differently than how they perceive how everyone else sees them. Minorities who are affected by imposter syndrome often fear that others will see them as undeserving of their position and presence which then can expose them as imposters, causing them to lose their status (Ewing et al. 2016). Feelings of being an imposter can result in diminished -self--efficacy, increased cognitive load, and decreased performance. Moreover, these feelings may lead many students and professionals to reconsider a career in S&T.

As authors of this paper, we come together with our shared, painful experiences acquired from participating in STEM professional organizations and conferences, to bring enlightenment to how these organizations can play a pivotal role in enhancing minority participation in STEM. Our collected experiences in two different fields has allowed us to synergize and gain a unique perspective that we will share with the hopes that the recommendations we make to professional societies will have a positive impact for all in every facet of the STEM professions.

## Where do professional societies fit in?

The truth of the matter is, we occasionally still encounter conference staff at majority meetings whose demeanors switch from cheerful and friendly to dry and deadpanned when we step forward to pay registration or to ask a simple question. Sadly, it is not limited to the race of the staffer. Many African American staffers mimic the worst of the behavior that we have seen. It may be worth noting that we consciously try to appear non-threatening at professional meetings by wearing professional attire including a sport jacket or suit, be reasonably clean-cut, and utilize proper grammar without ebonics. Even though we approach others with a polite tone and demeanor, we elicit a reaction as if we do not belong. We do not believe that these types of responses are personally motivated. Rather, these responses emanate from an unconscious confirmation bias, which is where one seeks to confirm one's beliefs or hypotheses, while giving disproportionately less consideration to alternative possibilities. These are common responses towards minority students and professionals from individuals who have had limited exposure or interests in interacting with people from outside their comfort zones. We also believe that these can be serious demotivating factors when it comes to minorities making choices about joining and participating in professional organizations.

In our current "post-racial" America, younger generations have grown up with a black President of the United States and what appears to be a plethora of equal opportunities and racial equity. The belief in a meritocracy, that if we perform the same, then the reward will also be the same takes precedence in their outlook. Many students and young professionals have not experienced the regularity of visceral racism of past decades nor are they as sensitive to the systemic and neoliberal racism of the current, post-racial era. This may inhibit their ability to navigate the more nuanced racial terrain. That occurs, when they enter a professional spaces in which they suddenly feel alienated, ignored, or subject to continuous micro-aggressions, they may fall away because of the hurtful, negative interactions and experiences. By maintaining a culture where this sort of treatment is normalized, we may be contributing to the attrition of the student and professional participation in STEM.

To persist as a career minority scientist, one often becomes acclimated to being the only black person in a room. In retrospect, this acclimation includes being the only minority in a room in addition to being the only minority in a hostile room. Our sense is that today's students coming from historically black colleges, minority serving institutions, or predominantly white institutions may never have experienced the kinds of overt racism that we grew up with in previous eras, nor the racism in its plethora of its distinguishable forms that developed as respectability politics became more of the norm. Subtle micro-aggressions can make professional meetings so overwrought with an inclement environment that many do not consider joining or attending on a regular basis. In our generation, we were told "to run faster, jump higher, and work longer and harder to gain the 'same' accolades as our white counterparts". To maintain this level of activity, many often complain of psychological burnout. We have continued to attend and participate in a variety of professional meetings and, unfortunately, the acceptance of diverse participants remains an issue that has yet to be resolved.

Attending a large professional meeting for the first time as a student can be quite daunting for anyone without adequate social networks or the natural ability to negotiate large groups of professionals angling for career advancement, renewing old acquaintances, and promoting their agendas. Establishing professional connections faced by a young scientist can be especially challenging for a lone African American, Latino/a, Native American, or person from another minority group. This has nothing to do with self-confidence or competence. It has far more to do with how misperceptions and discomfort with interacting with minority students or professionals can lead to unwelcoming and alienating experiences.

**Three common myths**

The three myths that we have selected for this article are by no means the only misperceptions limiting progress. These 3 are selected based on their prominence in the fields, which we are familiar.

*Myth 1*: There is no problem with underrepresentation. The degree of minority representation in the field is consistent with the existing pool of members of that group who have the talent, interest, and/or desire as attracted by market needs.

*Myth 2:* The "white savior model" is the most effective pathway towards success. In other words, majority entities can solve the diversity problem with their superior intellectual capital, post-racial politics, and research on minority issues sans people of color.

*Myth 3:* Diversity is equivalent to having more white females in the organization. Gender balance in a white monoculture simply preserves white privilege and exacerbates racial division and distrust in the long run as discussed in the rich literature of Bell Hooks and Kimberlé Williams Crenshaw (Hooks 1984; Hooks 2000; Crenshaw 2016).

There may be more myths, beliefs, and personal experiences that influence our perceptions and our actions. However, these three myths reflect a broad swath of the perceptual limitations that frame the responses to calls for increased diversity in STEM. Below, we will attempt to undertake a holistic approach to offering suggestions on how these myths can be overturned and a whole new approach to inclusion focusing on professional societies might have a positive impact on the trends discussed in the introduction. Going forward, we need to be aware of our beliefs as well as the beliefs of others so that we can all come together to rally for enhanced inclusivity in STEM. We now outline a few of the barriers that may be used to support and entrench the aforementioned myths in our professional culture. These barriers may represent the support scaffolding for the confirmation bias present in our professional societies. By better understanding these barriers and their underlying subtleties, perhaps we can all work towards making robust strides toward achieving equity in STEM.

**Ten potential barriers linked to the myths**

1. **Most people want to avoid conversations about race like the plague.** . The tendency to associate racial stereotypes with overt racism (now regarded as a fringe response to diversity in our "post-racial society" despite recent events) is a common defense mechanism among non-blacks. Even well-meaning people who often populate diversity committees may be uncomfortable in directly addressing racial issues within the organization. This becomes especially difficult when those exhibiting the biases or denying that structural biases exist are viewed as powerful figures in the organization.
2. **Difficult conversations about racial diversity often get diverted into conversations about gender diversity.** There are a variety of reasons for this diversion including the fact that diversity can mean a lot of different things to different subgroups/people. This can lead to discussions of alternate types of "isms" and false equivalences (classism, sexism, genderisms, etc.) – all of which are important but shifting focus rarely leads to a solved problem and completely fails to address ongoing biases. For example, stereotypes often present white females as less threatening, more approachable, and possessing a "shared" experience – especially in a racially homogenous context. This does not extend to women of color and men of color, especially those who fail to assimilate to the cultural norms of the society.

3. **The unspoken paradoxes associated with diversity events at professional meetings.** One paradox is that when one calls an event a "diversity event", then members of the monoculture (and even some from subgroups) may perceive that the event is exclusive and not inclusive to the entire community. Others may be hesitant to enter the space because they interpret diversity as being the "non-white" event and that they might be unwelcome interlopers. As a consequence, the event can take on the appearance of actually being an exclusive event.
4. **Many successful diversity events are attended because they are popular and hence, may attract individuals who do not sincerely support diversity efforts.** As a consequence, the initial motivating focus of the diversity efforts may be derailed and morph into a mission inconsistent with the spirit of the event. This is an example of virtue-signaling which can be defined as the conspicuous expression of inflated ethical values by an individual performed with the intent of presenting themselves as more enlightened than the broader social group to which they belong.
5. **Basic misconceptions about diversity in an organization**. Many people still do not understand that diversity is not solely about fairness and overcoming personal biases. It is more about benefitting from the full potential of all skilled professionals and becoming a better organization. We must remind ourselves that it is grounded in professional and scientific ethics. Diversity is about making sure that our occupational landscape is inclusive to all persons and ensures opportunities for students from all backgrounds in preparation for careers that may not yet exist.
6. **Deconvolving scientific objectivity: Racism and the post-racial stance.** Many scientists automatically (dogmatically) assume that they are unbiased, objective, and colorblind. That is, they believe that the very nature of their profession may place them beyond decisions based on racial stereotypes. Unfortunately, our professional arena can be rife with intelligent people in various stages of denial and adamant that no such biases exist. As scientists, we often tend to "dehumanize" ourselves while elevating the science. However, our human tendencies tend to dictate interpersonal interactions that are not included in the scientific method and can impact our judgment when we evaluate the quality of scientific works of others.
7. **Organizations do not naturally diversify.** In stark contrast, as humans, we tend to bend towards nepotism, professional incest, and traditional social networks. Minority students and newcomers can be preferentially excluded when these practices are entrenched. Racial amnesia can color conversations concerning diversity in professional spaces and can lead to the creation of blind spots. We share many ugly realities such as many of those in leadership positions gained their education during the "segregated sixties" or before. Ignoring or choosing not to confront past discrimination inhibits honest, but uncomfortable discussions about systemic biases. Colleagues have related instances during which individuals in leadership positions have intimated a nostalgia for the days when "the troublesome whiners" were not asking for handouts. Actually, nobody is asking for a handout. Diversity efforts may get channeled through standing committees of the professional organization. If these committees have other priorities, they will inevitably seek to either continue their standard modus operandi or argue pragmatism as a reason to fail to pursue diversity issues.
8. **Professional organizations often ask the wrong questions about improving diversity.** Often we get asked, "Why don't they [minorities] come?" We often get asked this when there have been historically low participation from people of color in scientific events. Our response, "What can you do to modify to make the event more enticing for them to come? How can we make them feel as welcome as everyone else? How can we demonstrate their tangible benefits for joining?"
9. **Change is hard work.** If what is being done now is not achieving the desired diversity, perhaps something different should be done. Addressing the resistance to change with the direct implications of how a lack of change can continue the plight of having barriers creep into our spaces will hopefully provide sufficient leverage. One has to work hard and consciously (strategically) to ensure that diversity events do not get labeled or perceived as exclusive and, more importantly, meet specific goals.

10. **Racial diversity is not a priority issue of all non-Caucasian professionals.** In fact, many "minorities" are repelled by the notion of diversity because of their own misleading about meritocracy and even basic misunderstandings about diversity. It is important to have people of all backgrounds sharing a common goal and vision of what diversity for the organization means in a practical sense, i.e. beyond a standard diversity statement.

**Potential strategies to enhance diversity in professional organizations**

One of the current efforts to address the diversity problem in Geoscience organizations evolved into the Colour of Weather Receptions at the American Meteorological Society Annual (AMS) meetings (Joseph et al. 2008). For the past twenty years, an informal networking reception has been hosted at this meeting focused on making it a more receptive and comfortable space for minority students and young professionals. These receptions grew out of ad hoc gatherings and get-togethers of no more than a handful of like-minded colleagues and grew to incorporate students of color from small meteorology programs at HBCUs and gradually expanded as networks and programs at minority-serving institutions continued to develop. The message at these receptions was never isolation. In fact, it openly encouraged and solicited leadership from the AMS to attend as a critical part of the inclusion process. Over the years, the event grew from ten or fifteen individuals to attendance exceeding two hundred and fifty participants. Colour of Weather has always sought to provide a safe space for provocative, often difficult, but necessary conversations about stereotypes and stereotype threats, overcoming systemic barriers to success, and for making personal contacts that extend beyond typical comfort zones. The two-hour receptions are now structured such that students can maximize professional interactions while engaging decision-makers (leadership from the AMS, federal agencies, private sector, and academic sector) in ice-breakers and activities that address bias and issues of diversity and inclusion as well as areas of common interest in atmospheric sciences, internship and training opportunities, and professional mentoring. The goal is for every minority student in attendance to leave the meeting with a stronger professional network than they possessed when they arrived. It is also designed to lower the activation barrier limiting interactions between students of color and the professional membership. We have been fortunate that the reception has become a regular part of the AMS meeting schedule and one of the highly sought out events of the opening day. Our support and attendance from AMS leadership has been a critical part of the success as has been the commitment of NASA and NOAA leadership. The regular attendance of these representatives have lent credence to the mission of the reception and has inspired even greater attendance from individuals who might have blown off the meeting otherwise.

In early 2017, the AMS President-Elect, Matthew Parker, approached one of the authors with openness and genuine interest to find what is sorely needed to enhance minority participation at AMS. As the Colour of Weather program was drawing to a close, he asked Morris, "How do we catalyze this energy and inclusiveness so that it expands throughout the society? We need these young people precisely because they are the future of this organization. I don't know how we do it but I'd like to work with you to realize this change." The openness and honesty that Parker conveyed in admitting that he did not have the answers but he actively sought advice and input was highly commendable. This showed true leadership of being able to identify potential cultural barriers and schisms, and seek the input of key members to find solutions so that the organization can thrive for all. By taking on the responsibility of making a change personally and with sincere honesty, solutions can grow in this fertilized land of strategic collaborations.

Given the impending demographic shifts in our Nation, it is in the best interests of all S&T professional societies to decrease the barriers that limit access and full inclusion in STEM. Professional organizations rely on membership. If the demographics of the membership of the professional societies does not evolve with the demographics of our society, then the sustainability of the professional societies may be jeopardized. Leadership in professional societies must also recognize that their voices may have an outsized impact on the culture of S&T workforce. Thus, by taking a definitive and proactive stance they can serve as positive change agents in society. Even though professional societies will not alone be able to completely resolve issues of racial equity in STEM, they may still be able to have a critical role in the solution.

Below, we distill several potential strategies that professional organizations can undertake in efforts to be part of the solution in improving diversity and inclusion in STEM and the S&T workforce.

1. **Racism was not created and cannot be solved by its victims.** Champions who represent the majority population and have significant stature are encouraged to be vocal about inclusion.
2. **Racism in society permeates the scientific fields and professional organizations.** Professional organizations may wish to take a proactive stance in trying to root out both active and passive racism in their meetings and organizational policies and practices. This may include nominations for awardees, memberships on boards, and key committee memberships.
3. **A possible key to fostering inclusion and diversity is helping minority students to establish robust social networks.** These network allows newcomers of all backgrounds to learn how to navigate in what can often be "hostile" environments for individuals who do not conform to perceived status quo.
4. **The purpose of diversity events should be recognized and may necessitate invitations to facilitate inclusiveness.** Further, it would be better if successful diversity events could be earmarked by tangible evidence that the organizational culture is changing for the better.
5. **The voice that advocated diversity should come from those who have earned a respected professional reputation.** This reputation can be achieved through tradition markers of academic success so that the voice can overcome skepticism of professional peers and directly address misgivings regarding diversity and exceptionalism. Often, given the current racial climate which can cloud our assessment of others, this voice should come from either a white or Asian male. This is not to say that more voices can bring harmony to STEM, for all are needed to push the infrastructure to be more inclusive.
6. **True partnerships between HBCUs as well as African American (e.g. Association of Black Geoscientists), Hispanic, and Native American scientific organizations and larger professional organizations should be forged to enhance participation, inclusion, and recognition.** An example could be reciprocal arrangements involving reduced membership fees for those belonging to minority professional societies. Minority scientists are often financially-stretched by attempting to maintain memberships in multiple organizations. At the same time, they are challenged to maintain time-commitments to committees in both organizations, while maintaining scientific productivity.
7. **Minority students and professionals should be knowledgeable on how to evaluate opportunities.** Their naive enthusiasm could be exploited by diversity programs seeking to fulfill a diversity requirement.
8. **Difficult conversations about the culture and environment of organizations must be undertaken and sustained.** The issue of diversity has been around for some time and as such, will take time and effort to be resolved.
9. **Access to leadership and leadership opportunities is a critical part of the solution.** All need to know and feel that they are a sincere member of the organization as this builds trust among the members that can help propagate the scientific community forward.
10. **Professional societies should not assume what people of color need.** Conversations and information-gathering should be undertaken with various groups or individuals. Finding these groups or individuals may take time and there may need to be some trust developed so that a sincere discussions can ensue. Note that this may require a degree of humility and openness to move forward on solutions that each party may not fully understand.

## Conclusion

We do not claim to have all of the answers in this paper. Yet we seek to inspire conversations that can synergize the strengths that spawn from our diversity which broaden our scientific communities. We conclude with the words of Lee Lorch (Case 1996), who continues to speak to all of our scientific societies' conscience.

> "Do we care only for those already in the profession? Is the Society willing to accept the present nearly complete exclusion from our mathematical manpower pool of Black America?
>
> Mathematicians and non-mathematicians alike await the answers to these questions. The Council [of all scientific societies] has the task of providing leadership."

## References

1. American Geophysics Institute. (2016) Status of the Geoscience Workforce 2016 https://www.americangeosciences.org/workforce/currents/underrepresented-minorities-us-workplace (accessed 4 Sept. 2017)
2. Ashkenas, J. Park H., and Pearce A. (2017) Even With Affirmative Action, Blacks and Hispanics Are More Underrepresented at Top Colleges Than 35 Years Ago, *New York Times*. https://www.nytimes.com/interactive/2017/08/24/us/affirmative-action.html (accessed 4 Sept. 2017)
3. Case B.A., Ed. (1996) *A Century of Mathematical Meetings*. American Mathematical Society, Providence, RI.
4. Crenshaw K. (2016) The Urgency of Intersectionality TED talk. https://www.ted.com/talks/kimberle_crenshaw_the_urgency_of_intersectionality (accessed 4 Sept. 2017)
5. Ewing K. M., Richardson T. Q., James-Myers L., and Russell R. K. (2016) The Relationship between Racial Identity Attitudes, Worldview, and African American Graduate Students' Experience of the Imposter Phenomenon, *Journal of Black Psychology,* **22**(1), 53-66.
6. Fiske S. T., Cuddy A. J. C., Glick P., and Xu J. (2002). A model of (often mixed) stereotype content: Competence and warmth respectively follow from perceived status and competition. *Journal of Personality and Social Psychology,* **82**(6), 878-902.
7. Goldberg T. (2015) *Are We All Post-Racial Yet?* John Wiley & Sons, New York, NY.
8. Hooks B. (1984) *Feminist Theory: From Margin to Center,* South End Press, Cambridge, MA.
9. Hooks B. (2000) *Feminism is for Everybody: Passionate Politics,* South End Press, Cambridge, MA.
10. Joseph E., Morris V. and Glakpe E. (2008) The NCAS Colour of Weather Event. *BAMS* **89**(7), 1042-1043.
11. Liebermann, R., Ehm L., and Gwanmesia G.. (2016). Creating career paths for African-American students in geosciences. *Eos,* **97**, doi:10.1029/2016EO052099.
12. National Science Foundation. (2017) Women, Minorities, and Persons with Disabilities in Science and Engineering 2017. https://www.nsf.gov/statistics/2017/nsf17310/static/downloads/nsf17310-digest.pdf (accessed 4 Sept. 2017)
13. Pew Research Organization. (2017) A look at historically black colleges and universities as Howard turns 150. http://www.pewresearch.org/fact-tank/2017/02/28/a-look-at-historically-black-colleges-and-universities-as-howard-turns-150/ (accessed 4 Sept. 2017)
14. U.S. Bureau of Labor Statistics. (2017) STEM Occupations Past, Present, and Future. https://www.bls.gov/spotlight/2017/science-technology-engineering-and-mathematics-stem-occupations-past-present-and-future/pdf/science-technology-engineering-and-mathematics-stem-occupations-past-present-and-future.pdf (accessed 4 Sept. 2017)
15. U.S. Census Bureau (2015) https://www.census.gov/content/dam/Census/library/publications/2015/demo/p25-1143.pdf (accessed 24 July 2017)